\documentclass[11pt]{article}
\usepackage{amsfonts}
\usepackage{}
\usepackage{mathrsfs}
\usepackage{amsmath}
\usepackage{verbatim}

\usepackage{amsmath, amsthm, amssymb}
\usepackage{enumerate}
\usepackage{graphicx}
\usepackage{float}
\usepackage{hyperref}
\usepackage[margin=1in]{geometry}
\numberwithin{equation}{section}

\begin{document}
\newtheorem{theorem}{Theorem}
\newtheorem{proposition}{Proposition}
\newtheorem{lemma}{Lemma}
\newtheorem{cor}{Corollary}
\newtheorem{definition}{Definition}
\theoremstyle{definition}
\newtheorem{exmp}{Example}[section]
\newtheorem{remark}[theorem]{Remark}

\title{An Extended Discrete Hardy-Littlewood-Sobolev Inequality}
\author{Ze Cheng \\ 
Congming Li\footnote{
Research partially supported by NSFC-11271166, NSF-DMS-0908097, and NSF-EAR-0934647.} \\
Department of Applied Mathematics, \\
University of Colorado Boulder, CO 80309, USA }
\date{}

\maketitle

\begin{abstract}
Hardy-Littlewood-Sobolev (HLS) Inequality fails in the ``critical'' case: $\mu=n$. However, for discrete HLS, we can derive a finite form of HLS inequality with logarithm correction for a critical case: $\mu=n$ and $p=q$, by limiting the inequality on a finite domain. The best constant in the inequality and its corresponding solution, the optimizer, are studied. First, we obtain a sharp estimate for the best constant. Then for the optimizer, we prove the uniqueness and a symmetry property. This is achieved by proving that the corresponding Euler-Lagrange equation has a unique nontrivial nonnegative critical point. Also, by using a discrete version of maximum principle, we prove certain monotonicity of this optimizer.
\end{abstract}

\section{Introduction}
The well-known Hardy-Littlewood-Sobolev (HLS) inequality states that
\begin{equation}\label{HLS}
\int_{\mathbb R^n} \int_{\mathbb R^n} \frac{f(x)g(y)}{|x-y|^{\mu}}\,dx \,dy \leq C_{p,\mu,n}\|f\|_p\|g\|_q
\end{equation}
for any $f \in L^p(\mathbb{R}^n)$ and $g \in L^q(\mathbb{R}^n)$ provided that
\[ 0 < \mu < n, 1 < p,q < \infty \text{ with } \frac{1}{p} + \frac{1}{q} + \frac{\mu}{n} = 2. \]

$C_{p,\mu,n}$ is the best constant for \eqref{HLS}, and proved by Lieb \cite{Lieb83} that, such $C_{p,\mu,n}$ and corresponding maximizing pair $(f,g)$ exists. In particular, Lieb also gave the explicit $f$ abd $C_{p,\mu,n}$ in the case $p=q$. The method Lieb used was to examine the Euler-Lagrange equation that the maximizing pair $(f,g)$ satisfies with some techniques to exploit the symmetry of $f$. This idea is inherited in \cite{LV11} and here to find the sharp estimate of best constant of a finite form of HLS in a critical case: $p=q=2$, and hence $\mu=n$.

Following the idea that the maximizer of HLS satisfies corresponding E-L equations, the study of the HLS inequality and weighted inequality later generalized by Stein and Weiss \cite{SW58} is naturally related to the studies of various of integral equations. For recent results, see \cite{CL91, CLO06, CL08, CL05} and a brief summary can be found in \cite{CJLL05}. These works have studied regularity and radial symmetry of solutions of such integral systems, and introduced a method of moving plane in an integral form which is proved to be a powerful tool. In \cite{Hang}, the result of integral system corresponding to HLS \eqref{HLS} is improved to all cases, i.e. the condition $p,q \geq 1$ is removed. In this paper, we do not use the method of moving plane directly, but borrowing its idea, we use a maximum principle to deal with a discrete problem and prove the symmetry of the solution.

First, let's have a look at the discrete and 1-dimensional version of HLS inequality \eqref{HLS}, the Hardy-Littlewood-P\'{o}lya (HLP) Inequality \cite{Hardy:52a}: if $ a \in l^p(\mathbb{Z})$ and $b \in l^q(\mathbb{Z})$ and
\begin{equation*}
0 < \mu < 1, 1 < p,q < \infty \text{ with } \ \frac{1}{p} + \frac{1}{q} + \mu = 2,
\end{equation*}
then
\begin{equation}\label{HLP}
\displaystyle \sum_{r\neq s} \frac{a_r b_s}{|r-s|^{\mu}} \leq C\|a\|_{p}\|b\|_{q}
\end{equation}
where $r,s\in \mathbb{Z}$ and the constant $C$ depends on $p$ and $q$ only.

For this HLP inequality \eqref{HLP}, let's consider the critical case: $p = q = 2$ and $\mu = 2 - \frac{1}{p} - \frac{1}{q} = 1$, for which the original HLP fails, but we can compromise and get a finite form of HLP. In \cite{LV11}, the inequality is extended to the critical case as: If $a, b \in l^p(\mathbb{Z})$, then
\begin{equation}\label{HLPextended1d}
\sum_{r\neq s,1\leq r,s \leq N} \frac{a_r b_s}{|r-s|} \leq \lambda_N\|a\|_{2}\|b\|_{2}.
\end{equation}
where $\lambda_N$ is the best constant for \eqref{HLPextended1d}, and $\lambda_N=2\ln N + O(1)$.

\begin{remark}
 One of the reasons that we consider discrete version of HLS instead of the original inequality is, when $\mu=1$ the integrand on the left side of HLS \eqref{HLS} is not always integrable on a finite domain for $L^p$ functions. So it is not as convenient to extend 1-dimension HLS inequality \eqref{HLS} to the critical case in a similar finite form as to extend HLP \eqref{HLP} to \eqref{HLPextended1d}.
\end{remark}

As for the high dimensional discrete HLS, if $a, b \in l^{p}(\mathbb{Z}^n)$, and
\begin{equation*}
0 < \mu < n, 1 < p,q < \infty \text{ with } \ \frac{1}{p} + \frac{1}{q} + \frac{\mu}{n} = 2,
\end{equation*}
then
\begin{equation}\label{discreteHLS}
\displaystyle \sum_{r\neq s} \frac{a_r b_s}{|r-s|^{\mu}} \leq C\|a\|_{p}\|b\|_{q}
\end{equation}
where $r,s\in \mathbb{R}^n$ and the constant $C$ depends on $p$ and $q$ only. We can extend \eqref{discreteHLS} to a finite form in the corresponding critical case: $p = q = 2$ and $\mu = n$, in the following way:
\begin{theorem}\label{BCEstimate}
If $r,s\in\mathbb{R}^n$ and $1\leq r_i, s_i\leq N$ where $r_i,s_i$ are integers and $1\leq i\leq n$, then $a_r,b_s\in \mathbb{R}^L$, where $L=N^n$. let
\begin{equation}\label{BC}
    \lambda_N=\max_{\|a\|_2=\|b\|_2=1}\sum_{r\neq s}{\frac{a_rb_s}{|r-s|^n}}
\end{equation}
So, we have an extension of HLS inequality
\begin{equation}\label{HLPextended}
    \sum_{r\neq s}{\frac{a_rb_s}{|r-s|^n}} \leq \lambda_N\|a_r\|_2\|b_s\|_2
\end{equation}
where the two statements below holds
\begin{description}
  \item[(i)]  $|S^{n-1}|\ln N-o(\ln N)<\lambda_N < |S^{n-1}|\ln N+o(\ln N)$.
  \item[(ii)] $\exists! \overline{a^N}=\overline{b^N}$ and $\|\overline{a^N}\|_2=1$ such that the equality in \eqref{HLPextended} holds, and $\overline{a^N}\in \mathbb{R}^L_+$ where $L=N^n$.
\end{description}
\end{theorem}

Let's call the triplet $(\overline{a^N}=\overline{b^N},\lambda_N)$ the optimizer of \eqref{HLPextended} since it is unique, and there are some properties of the optimizer. First, as a consequence of the uniqueness, we have symmetry property of the optimizer in the following sense,
\begin{theorem}\label{symmetry}
Let $(\overline{a^N},\lambda_N)$ be the optimizer. $\Phi:S\rightarrow S$ is an isometric map, where $S=\{r\in\mathbb{R}_+^n|1\leq r_i\leq N\}$. Then $\overline{a^N_{\Phi(r)}}=\overline{a^N_r}$.
\end{theorem}

Second, the optimizer has certain monotone decaying property. For convenience of writing, let's change the range of $r_i$ from $[1,N]$ to $[-N,N]$, which makes no essential change to the results above, and we have the monotone decaying property for this special case,
\begin{theorem}\label{decay}
If $(\overline{a^N},\lambda_N)$ is the optimizer and $r\in\mathbb{R}^n$, $-N\leq r_i\leq N$ for $1\leq i \leq n$, then $a\in\mathbb{R}_+^L$, where $L=(2N+1)^n$, and $\overline{a^N}$ has a monotone decaying property from its central element: For $1\leq i \leq n$,
\begin{equation}\label{DecayingProperty}
    \left\{ \begin{aligned}
        \overline{a^N_{(r_i,r')}} &\leq \overline{a^N_{(r_i-1,r')}},  1\leq r_i\leq N \\
        \overline{a^N_{(r_i,r')}} &\geq \overline{a^N_{(r_i-1,r')}},  -N+1\leq r_i\leq 0
        \end{aligned}
\right.
\end{equation}
\end{theorem}

To prove theorem \ref{decay}, we use the following maximum principle,
\begin{theorem}[Maximum Principle]\label{MP}
Let $\mathbb{R}_+^L$ be the positive cone in $\mathbb{R}^L$, i.e., if $a \in \mathbb{R}_+^L$ then every element of $a$ is positive. Suppose a linear equation:
\begin{equation}\label{LinearEq}
    u = Au+f
\end{equation}
where $A:\overline{\mathbb{R}_+^L}\rightarrow\overline{\mathbb{R}_+^L}$ with $\|A\|_2<1$, and $f\in\overline{\mathbb{R}_+^L}$, then $\exists ! u$ satisfies \eqref{LinearEq} and $u \in \overline{\mathbb{R}_+^L}$. In other words, $(I-A)^{-1}\in\overline{\mathbb{R}_+^{L\times L}}$.
\end{theorem}

This Maximum Principle follows directly from standard contracting mapping iteration. It is a discrete version of maximum principle analogous to the usual versions in PDE. To see this, let's look at a typical maximum principle: let $\Omega\subset\mathbb{R}$ be an open bounded and connected domain with smooth boundary $\partial\Omega$. Let $u \in C^2(\Omega)\cap C(\overline \Omega)$ be a solution of following equation,
\begin{equation}\label{EllipticEq}
    \left\{ \begin{aligned}
        -\Delta u &= f \geq 0 \ \text{in} \ \Omega \\
        u &= 0 \ \text{on} \ \partial\Omega
        \end{aligned}
\right.
\end{equation}
Then by maximum principle $u \geq 0$ in $\Omega$. Actually, by strong maximum principle, $u>0$ or $u\equiv0$ in $\Omega$.

So, theorem \ref{MP} is indeed saying that if $(I-A)u=f \in \overline{\mathbb{R}_+^L}$, then $u \in \overline{\mathbb{R}_+^L}$. Corresponding to strong maximum principle, in theorem \ref{MP} if every entry of $A$ is strictly positive, it is easy to see that $u \in \mathbb{R}_+^L$. For more general symmetric linear operators, there is also maximum principle, and one can check \cite{Kigami01} for details.

\section{Best Constant Estimate in High Dimension Space}

{\bf{Proof of part (i) of theorem \ref{BCEstimate}.}} {\bf{Step 1.}} $\lambda_N\geq |S^{n-1}|\ln N-o(\ln N)$.

Let $a=b$, and
\begin{equation}
    a_r = N^{-\frac{n}{2}},1\leq r_i\leq N, 1\leq i\leq n
\end{equation}
So, $\|a\|_2=1$.

By the definition of $\lambda_N$, we have
\begin{align*}
    \lambda_N &\geq \sum_{r\neq s}{\frac{a_ra_s}{|s-r|^n}}=N^{-n}\sum_{r\neq s}{\frac{1}{|s-r|^n}} \\
            & =N^{-n}\{2^n\sum_{x_n=1}^{N-1}\cdots\sum_{x_1=1}^{N-1}\frac{(N-x_1)\cdots (N-x_1)}{(x_1^2+\cdots x_n^2)^{\frac{n}{2}}}-o(N^n\ln N)\} \\
            & \geq (\frac{N}{2})^n\int_0^{\frac{\pi}{2}}\cdots\int_0^{\frac{\pi}{2}}\int_1^N \frac{(N-r\cos\phi_1)\cdots}{r^n}r^{n-1}drd\phi_1\cdots d\phi_n -o(\ln N) \\
            & = (\frac{2}{N})^n|S^{n-1}|2^{-n}N^n\ln N-o(\ln N) \\
            & = |S^{n-1}|\ln N-o(\ln N)
\end{align*}

{\bf{Step 2.}} $\lambda_N\leq |S^{n-1}|\ln N+o(\ln N)$ \\
 Let $J(a,b)=\sum_{r\neq s}\frac{a_rb_s}{|r-s|^n}$. Hence, $\lambda_N=\max_{\|a\|_2=\|b\|_2=1} J(a,b)$, i.e. we will maximize $J(a,b)$ under the constraints $\|a\|_2=\|b\|_2=1$ (in fact, we use $\frac{1}{2}\|a\|_2^2=\frac{1}{2}\|b\|_2^2=\frac{1}{2}$).  Therefore, we conduct Euler-Lagrange equations and by compactness: $\exists \|\overline{a^N}\|_2=\|\overline{b^N}\|_2=1$ such that $\lambda_N=J(\overline{a^N},\overline{b^N})$ and,
\begin{equation}\label{EulerLagrangeOringin}
    \left\{ \begin{aligned}
         \lambda_1 \overline{a^N_r} &= \sum_{s\neq r}\frac{\overline{b^N_s}}{|s-r|^n} \\
                  \lambda_2 \overline{b^N_s} &= \sum_{r\neq s}\frac{\overline{a^N_r}}{|r-s|^n}
                          \end{aligned} \right.
\end{equation}
where $r,s\in\mathbb{R}^n$ and $1\leq r_i, s_i\leq N$.

For convenience, write \eqref{EulerLagrangeOringin} in matrix form,
\begin{equation} \label{EulerLagrangeMatrix}
\left\{ \begin{aligned}
         \lambda_1 \overline{a^N} &= A\overline{b^N} \\
                  \lambda_2 \overline{b^N} &= A\overline{a^N}
                          \end{aligned} \right.
\end{equation}
Left multiply the first equation of \eqref{EulerLagrangeMatrix} by $a^T$, the second equation by $b^T$, and by the fact that $A$ is symmetric and $\|\overline{a^N}\|_2=\|\overline{b^N}\|_2=1$, one sees that
\begin{align*}
    \lambda_1 &= \lambda_1\|\overline{a^N}\|_2^2 = \overline{a^N}^TA\overline{b^N}= J(\overline{a^N},\overline{b^N}) \\
            &=\overline{b^N}^TA^T\overline{a^N} = \lambda_2\|\overline{b^N}\|_2^2 =\lambda_2
\end{align*}
and since $\lambda_N=J(\overline{a^N},\overline{b^N})$, we have $\lambda_1=\lambda_2=\lambda_N$. \\
Now, let $b_{s_0}=\max{|\overline{a_r^N}|,|\overline{b_s^N}|}>0$, so, $\overline{b_{s_0}^N} \lambda_N=\sum_{r\neq s_0}\frac{\overline{a^N_r}}{|r-s_0|^n}$, which leads to
\begin{align*}
    \lambda_N &= \sum_{r\neq s_0}\frac{\overline{a_r^N}}{b_{s_0}|r-s_0|^n} \leq \sum_{r\neq s_0}\frac{1}{|r-s_0|^n} \\
            & \leq \sum_{r\neq (\frac{N}{2},\cdots,\frac{N}{2})=m_0}\frac{1}{|r-m_0|^n} \\
            & \leq \int_\Sigma\int_1^{\frac{\sqrt{2} N}{2}}\frac{1}{r^n}r^{n-1}drd\sigma \\
            & \leq |S^{n-1}|(\ln\frac{\sqrt{2} N}{2}) = |S^{n-1}|(\ln N +\frac{1}{2}\ln2) \\
            & = |S^{n-1}|\ln N + o(\ln N)
\end{align*}

Part (ii) will be shown later in section 3.
$\square$

\begin{lemma}\label{SingleSign}
If $(a,b,\lambda_N)$ satisfies $\|a\|_2=\|b\|_2=1$ and makes the equality of \eqref{HLPextended} hold, then $a,b\in\overline{\mathbb{R}_+^L}\cup\overline{\mathbb{R}_-^L}$.
\end{lemma}

Notice that if there is a sign change among the elements of $a$ and $b$, $(a,b)$ must not be an optimizer since $|\sum a_ib_i|<\sum |a_i||b_i|$. So the lemma holds, and it means that we can assume the triplet $(\overline{a^N},\overline{b^N}, \lambda_N)$ above to satisfy $\overline{a^N},\overline{b^N}\in\overline{\mathbb{R}_+^L}$.

Now, let's introduce a notation,
\begin{definition}\label{optimizer}
$(a,b, \lambda_N)$ such that
 \begin{itemize}
   \item $\|a\|_2=\|b\|_2=1$
   \item $a,b\in \overline{\mathbb{R}_+^L}$
   \item The equality of \eqref{HLPextended} holds
 \end{itemize}
is called an optimizer or solution of optimization of \eqref{HLPextended}.
\end{definition}
Obviously, $(\overline{a^N},\overline{b^N}, \lambda_N)$ is an optimizer. Next, we are going to prove part(ii) of theorem \ref{BCEstimate}, i.e., the optimizer is unique in positive cone and $\overline{a^N}=\overline{b^N}$.

\section{Uniqueness of The Optimizer}

From previous discussion we see that, an optimizer of \eqref{HLPextended}, $(\overline{a^N},\overline{b^N}, \lambda_N)$, satisfies Euler-Lagrange equations\eqref{EulerLagrangeOringin}. We are going to show the optimizer is unique in positive cone by showing the solution of the Euler-Lagrange equations in the positive cone $\mathbb{R}_+^L$ where $L=N^n$ is unique. Considering the following equations,
\begin{equation} \label{EulerLagrange}
\left\{ \begin{aligned}
         \lambda_1 a_r &= \sum_{s\neq r}\frac{b_s}{|s-r|^n} \\
                  \lambda_2 b_s &= \sum_{r\neq s}\frac{a_r}{|r-s|^n}
                          \end{aligned} \right.
\end{equation}
 where $\|a\|_2=\|b\|_2=1$, $r=(r_i)\in\mathbb{R}^n$, and $1\leq r_i\leq N$, $1\leq i\leq n$. $a,b\in\mathbb{R}^L$, where $L=N^n$. By lemma \ref{SingleSign}, we only need to study solution of \eqref{EulerLagrange} in the positive cone $\overline{\mathbb{R}_+^L}$.

In the proof, we will use the following simple map,
\begin{definition}\label{absoluteMap}
Let $T:\overline{\mathbb{R}^L}\rightarrow \overline{\mathbb{R}_+^L}$ such that $(Ta)_i=|a_i|$ for $1\leq i \leq L$.
\end{definition}

\begin{theorem}\label{uniqueness}
If $(a,b,\lambda_1, \lambda_2)$ is a solution of \eqref{EulerLagrange}, where $a,b \in \overline{\mathbb{R}_+^L}$, then $\lambda_1=\lambda_2=\lambda_N$, and $a=b\in\mathbb{R}_+^L$ is unique.
\end{theorem}

Proof. {\bf{Step 1.}} $\lambda_1=\lambda_2$.

This is similar to step 2 of theorem \ref{BCEstimate}. So, let $\lambda=\lambda_1=\lambda_2$.

{\bf{Step 2.}} $a,b\in\mathbb{R}_+^L$.

Since
\begin{align*}
    \lambda a_r &= \sum_{s\neq r}\frac{b_s}{|s-r|^n} \\
                     &= \frac{1}{\lambda}\sum_t\sum_{s\neq r,t}(\frac{1}{|r-s|^n}\frac{1}{|t-s|^n})a_t \\
                     &= \frac{1}{\lambda}\sum_t C(r,t)a_t
\end{align*}
we have $\lambda^2a=Ca$, where $C=A^TA$ and $A$ is a symmetric matrix. So C is non-negative definite. Since $C(r,t)>0$, $a \in \overline{\mathbb{R}_+^L}$ and $a\neq 0$ for $\|a\|=1$, the last term above is strictly positive. Therefore, $a,b\in\mathbb{R}_+^L$.

Let $0\leq\mu_1\leq\mu_2\cdots\leq\mu_L$ be the eigenvalues of $C$. Then $\exists \overline{\xi_L}\in\mathbb{R}^L$, s.t. $C\overline{\xi_L}=\mu_L\overline{\xi_L}$, and $\|\overline{\xi_L}\|=1$, and $\overline{\xi_L}\notin\overline{\mathbb{R}^L_-}$. We can assume the last property because eigenvectors appear in pairs with opposite signs. Also, by theory of adjoint operators, $\mu_L=\sup_{\|\xi\|=1}<\xi,C\xi>=<\xi_L,C\xi_L>$.

{\bf{Step 3.}} $\exists\overline{\xi_L}\in \mathbb{R}_+^L, \|\overline{\xi_L}\|=1$, and $\mu_{L-1}<\mu_L$.

First, $\exists\overline{\xi_L}\in \overline{\mathbb{R}_+^L}$. If not, then $\overline{\xi_L}\notin \overline{\mathbb{R}_+^L}\cup\overline{\mathbb{R}_-^L}$.

Then we have
\begin{align}
         \mu_L &= \overline{\xi_L}^TC\overline{\xi_L} \\
                     &< (T\overline{\xi_L})^TC(T\overline{\xi_L}) \\
                     &\leq \max_{\|\xi\|=1} \xi^TC\xi = \mu_L
\end{align}
where $T$ is defined in definition\ref{absoluteMap}.
A contradiction. So, $\exists \overline{\xi_L}\in\overline{\mathbb{R}_+^L}$, and since $C\overline{\xi_L}=\mu_L\overline{\xi_L}$, $\overline{\xi_L}\in\mathbb{R}_+^L$.

The argument above also shows that $\mu_{L-1}<\mu_L$. If not, $\mu_{L-1}=\mu_L$, then by a similar argument as above $\exists \xi_{L-1}\in\mathbb{R}_+^L$, s.t. $C\xi_{L-1}=\mu_L\xi_{L-1}$, and moreover $\xi_{L-1}\bot\xi_L$ which is impossible.

{\bf{Step 4.}} $a=b=\overline{\xi_L}, \lambda=\lambda_N=\sqrt{\mu_L}$.

Considering $\lambda^2a=Ca$,
\begin{enumerate}
  \item If $\lambda^2\neq\mu_L$, then $a\bot\overline{\xi_L}$. Since $a\in\mathbb{R}_+^L$ by step 2, this is impossible. So, $\lambda^2=\mu_L$.
  \item Since $Ca=\mu_La, C\overline{\xi_L}=\mu_L\overline{\xi_L}$, and by the fact that $\mu_{L-1}<\mu_L$ and $\|a\|=\|\overline{\xi_L}\|=1$, $a=\overline{\xi_L}$. Similarly, $b=\overline{\xi_L}$.
  \item If $(\overline{a^N},\overline{b^N}, \lambda_N)$ is an optimizer of \eqref{HLPextended} in the positive cone, it is a solution of \eqref{EulerLagrange}. So, $a=\overline{a^N}=b=\overline{b^N}$, $\lambda_N^2=\lambda^2=\mu_L$, and $\lambda,\lambda_N>0$, so $\lambda=\lambda_N$. $\square$
\end{enumerate}

\noindent
{\bf{Proof of part (ii) of theorem \ref{BCEstimate}.}} The same as the 3rd argument of step 4 above, since an optimizer $(\overline{a^N},\overline{b^N}, \lambda_N)$ is a solution of \eqref{EulerLagrange}, part (ii) follows from theorem \ref{uniqueness}.

\begin{remark}
At the time of this writing, thanks to Professor Dongsheng Li of Jiaotong University in Xi'an, we find that uniqueness follows directly from Perron's theorem \cite{Perron1907}. So the proof above can be much simplified.
\end{remark}

\begin{cor}\label{numdamonotone}
$\lambda$ is increasing as $N$ increases.
\end{cor}
Proof. Let $\lambda_N$ and $A_N$ be a solution and coefficient matrix of \eqref{EulerLagrange}. So,
\begin{align*}
    \lambda_N &= \max_{\|\xi\|=1}\xi^TA_N\xi =\overline{\xi_N}^TA_N\overline{\xi_N} \\
            &= (\overline{\xi_N},0)^TA_{N+1}(\overline{\xi_N},0) \\
            &< \max_{\|\xi\|=1}\xi^TA_{N+1}\xi = \lambda_{N+1}
\end{align*}
where $(\overline{\xi_N},0)$ means $(\overline{\xi_N},0)\in\mathbb{R}^L$ and $L=(N+1)^n$, and arranging $\overline{\xi_N}$ to take the first $N^n$ entries and stuffing the rest with zeros. Then calculate in blocks of matrices. $\square$

\section{Symmetry of The Optimizer}
From section 2 we see the uniqueness of the optimizer of \eqref{HLPextended} in positive cone. So, from this point, if it is clear in context, we use $(a=b,\lambda)$ instead of $(\overline{a^N},\overline{b^N}, \lambda_N)$ when referring the optimizer of \eqref{HLPextended} for simplicity.

\noindent
{\bf{Proof of theorem \ref{symmetry}.}} From \eqref{EulerLagrange} we have
\begin{equation*}
    \lambda a_r = \sum_{s\neq r}\frac{a_s}{|s-r|^n}
\end{equation*}
then
\begin{align*}
    \lambda a_{\Phi(r)} = \sum_{s\neq \Phi(r)}\frac{a_s}{|s-\Phi(r)|^n} = \sum_{t\neq r}\frac{a_{\Phi(t)}}{|\Phi(t)-\Phi(r)|^n} = \sum_{t\neq r}\frac{a_{\Phi(t)}}{|t-r|^n}
\end{align*}
So, $\bar{a}=(a)_{\Phi(r)}$ is also a solution to \eqref{EulerLagrange}. Then, by uniqueness of the solution, $\bar{a}=a$. So, $a_{\Phi(r)}=a_r$. $\square$

\begin{exmp}\label{exampleSymmetry}
If $a$ is an optimizer, then $a_{(r_i,r')}=a_{(N-r_i+1,r')}$ for $1\leq i\leq N$.
\end{exmp}

\section{Monotone Property of The Optimizer}
For convenience of writing, we change the range of $r_i$'s from $1\leq r_i\leq N$ to $-N\leq r_i\leq N$ which makes no change to the results above essentially.

\noindent
\\
{\bf{Proof of Theorem \ref{decay}.}}
We are only going to show \eqref{DecayingProperty} is true for $i=1$ for simplicity. Consider $d_r^{(1)}=a_{(r_1-1,r')}-a_{(r_1,r')}$, where $1\leq r_1\leq N$ and $-N\leq r_i\leq N$, $2\leq i\leq n$. So $d^{(1)}\in\mathbb{R}^{N(2N+1)^{(n-1)}}$. Then by applying theorem \ref{symmetry}, we have
\begin{align*}
    d_r^{(1)} &= \frac{1}{\lambda}(\sum_{s\neq (r_1-1,r')}\frac{a_s}{|s-(r_1-1,r')|^n}-\sum_{s\neq (r_1,r')}\frac{a_s}{|s-(r_1,r')|^n}) \\
     &= \frac{1}{\lambda}(\sum_{\substack{t=(t_1,t')\neq (r_1,r'), \\-N+1\leq t_1\leq N+1}}\frac{a_{(t_1-1,t')}}{|t-(r_1,r')|^n}-\sum_{s\neq (r_1,r')}\frac{a_{(s_1,s')}}{|s-(r_1,r')|^n}) \\
        &= \frac{1}{\lambda}(\sum_{\substack{t=(t_1,t')\neq (r_1,r'), \\1\leq t_1\leq N}}\frac{d_t^{(1)}}{|t-(r_1,r')|^n} + \sum_{\substack{t=(-t_1+1,t')\neq (r_1,r'), \\1\leq t_1\leq N}}\frac{-d_t^{(1)}}{|(-t_1+1,t')-(r_1,r')|^n} \\
        &+\underbrace{\sum_{\substack{t=(N+1,t')\neq (r_1,r')}}\frac{a_{(N,t')}}{|t-(r_1,r')|^n} - \sum_{t=(-N,t')\neq (r_1,r')}\frac{a_{(-N,t')}}{|t-(r_1,r')|^n}}_{f(r)}) \\
    &= \frac{1}{\lambda}(\sum_{\substack{(t_1,t')\neq r\\ 1\leq t_1\leq N  }}(\frac{1}{|(t_1,t')-r|^n}-\frac{1}{|(-t_1+1,t')-r|^n})d_t^{(1)} + \frac{-d_r^{(1)}}{|2r_1-1|^n} + f(r))
\end{align*}

 Also by theorem \ref{symmetry}, $a_{(N,t')}=a_{(-N,t')}$, easily one sees that $f(r)\geq 0$.

So, for $1\leq r_1\leq N$
\begin{align*}
    (\lambda+\frac{1}{|2r_1-1|^n})d_r^{(1)} &= \sum_{\substack{(t_1,t')\neq r\\ 1\leq t_1\leq N}}(\frac{1}{|(t_1,t')-r|^n}-\frac{1}{|(-t_1+1,t')-r|^n})d_t^{(1)} + f(r) \\
\end{align*}
Write the above equations in matrix form,
\begin{equation}\label{keyMP}
    d^{(1)} = Ad^{(1)} + F
\end{equation}
where $(F)_r=\frac{1}{(\lambda+\frac{1}{|2r_1-1|^n})}f(r)$, and
\begin{equation*}
A(r,t)=\left\{\begin{aligned}
    &\frac{1}{(\lambda+\frac{1}{|2r_1-1|^n})}(\frac{1}{|(t_1,t')-r|^n}-\frac{1}{|(-t_1+1,t')-r|^n}), &r\neq t\\
    & 0, &r=t
        \end{aligned}\right.
\end{equation*}

It is easy to see that entries of $A$ and $F$ are non-negative. So, $A:\overline{\mathbb{R}_+^L}\rightarrow\overline{\mathbb{R}_+^L}$, where $L=N(2N+1)^{(n-1)}$, and $F\in \overline{\mathbb{R}_+^L}$. Therefore, provided $\|A\|<1$, then by Theorem \ref{MP} (Maximum Principle) we get $d^{(1)}\in\overline{\mathbb{R}_+^L}$, hence \eqref{DecayingProperty} is proved. So, the only thing left to prove is $\|A\|<1$.

Notice that if $C,D$ are symmetric matrices such that $C,D:\mathbb{R}_+^L\rightarrow\mathbb{R}_+^L$, for some positive integer $L$, then $\|C\|\leq\|C+D\|$, because
\begin{align*}
    \|C\|=\max_{\|\xi\|=1}\xi^TC\xi=\bar{\xi}^TC\bar{\xi}\leq \bar{\xi}^T(C+D)\bar{\xi}\leq\max_{\|\xi\|=1}\xi^T(C+D)\xi = \|C+D\|
\end{align*}
Let
\begin{equation*}
C(r,t)=\left\{\begin{aligned}
    &\frac{1}{(\lambda+\frac{1}{|2r_1-1|^n})}\frac{1}{|(t_1,t')-r|^n}, &r\neq t\\
    & 0, &r=t
        \end{aligned}\right.
\end{equation*}
and
\begin{equation*}
D(r,t)=\left\{\begin{aligned}
    &\frac{1}{(\lambda+\frac{1}{|2r_1-1|^n})}\frac{1}{|(-t_1+1,t')-r|^n}, &r\neq t\\
    & 0, &r=t
        \end{aligned}\right.
\end{equation*}
So,
\begin{align*}
    \|A\| &\leq \|A+D\| =\|C\| \\
            &\leq \frac{1}{\lambda+\delta(N)}\|A_N\|
\end{align*}
where $A_N$ is the matrix of \eqref{EulerLagrange} of the case that $-N\leq r_i\leq N$ and $1\leq i\leq n$, so $\|A_N\|=\lambda$. So, $\|A\|<1$. $\square$

\bibliography{20130104bib}{}
\bibliographystyle{plain}

\end{document}